\documentclass[12pt]{article}
\usepackage{graphicx} % Required for inserting images
\usepackage{amssymb}
\usepackage{amsfonts}
\usepackage{amsmath}
\usepackage{tikz}

\title{An elementary way to find a counterexample to the Jacobian Conjecture}
\author{Arno van den Essen}
\date{12 September 2026}

\begin{document}

\maketitle

\section{Introduction}
On July 20 I was pleasantly surprised by an email of one of my former students Jan Schoone, who sent me the three dimendional counterexample to the Jacobian Conjecture, found by Levent Alpöge [1] a few hours earlier. On the other hand it was not at all clear how the author had found this example, other than knowing that he had used AI, more precisely Anthropic's AI model Claude Fable 5, which was released on June 9, 2026. In the meantime several authors commented on this result and found new counterexamples [3],[4],[7] and [8].
Since I have been working on this conjecture since 1985 I wanted to understand how this example could have been found without using AI. which is the aim of this note. Our paper is more modest than the papers cited above. Our aim is just to show how a counterexample can be found with only very elementary mathematics. Before we embark on this journey let us recall some basic facts and notations concerning polynomial maps.

\section{Polynomial maps}

A map $F=(F_1,\ldots,F_n):\mathbb{C}^n\rightarrow \mathbb{C}^n$ is called a {\em polynomial map} if each component $F_i$ is a polynomial in the variables $x_1,\ldots,x_n$ over the complex numbers $\mathbb{C}$. Such a map is called {\em invertible} if there  exists a polynomial map $G=(G_1,\ldots,G_n)$ such that the composition 

$$F\circ G:=(F_1(G_1,\ldots,G_n),\ldots,F_n(G_1,\ldots,G_n))$$
 is the identity map i.e.
$F_i(G_1,\ldots,G_n)=x_i$ for all $i$. A component of an invertible polynomial map is called a {\em coordinate}. It easily follows from the chain rule that a necessary condition for $F$ to be invertible is that its {\em Jacobian determinant} $det JF$ is a non-zero constant, where $JF$ denotes the Jacobian matrix of $F$. This condition is called the {\em Jacobian condition}. Polynomial maps satisfying the Jacobian condition are called {\em Keller maps}. It was Keller who in 1939 (see [5]) posed the question if every Keller map is invertible. This question became later known under the name {\em Keller's problem} but is more often cited as the {\em Jacobian Conjcture}.
For more details we refer to [2].

In order to find counterexamples to this conjecture we must be able to construct Keller maps. So in three variables we need to find three polynomials $f,g,h\in\mathbb{C}[x,y,z]$ such that $detJ(f,g,h)$ is a non-zero constant. We call such a triple a {\em Jacobian triple}. By developing this determinant according to the first row of the Jacobian matrix this implies in particular that the gradient of $f$ i.e. the row $(f_x,f_y,f_z)$ is a unimodular row. So to start finding a Jacobian triple one has to look at polynomials $f$ such that their gradient is  a unimodular row. Furthermore, to be sure that the triple  $(f,g,h)$ we construct is a counterexample we want to impose on $f$ an extra condition which guarantees that $(f,g,h)$ is not invertible. 

So what can we say about $f$ if $F=(f,g,h)$ is invertible? Well, composing $F$ with its polynomial inverse $G$ we deduce that $f(G_1,G_2,G_2)=x$. From this it follows easily that $f$ is irreducible. So if we can find a {\em reducible} polynomial $f$, having a unimodular row, which can be extended to a Jacobian triple $(f,g,h)$ we are sure to have found a counterexample to the Jacobian Conjecture. 

When on July 20 I sent Alpöge's counterexample to several of my former students, it was Stefan Maubach who phoned me and during our conversation made the remark that the third component $F_3$ of Alpöge's map is a so-called {\em fake coordinate}, i.e. a component of a Keller map that is not a coordinate. This is obvious since $F_3=x(2-3xy-x^2z)$ is reducible and Alpöge's map is a Keller map. The word {\em fake polynomial} reminded me of the polynomial $f=x+x^2y$ which is the simplest reducible polynomial such that its gradient $(f_x,f_y)$ is unimodular, in other words such that this row can be extended to an invertible $2\times 2$ matrix. If this extension could be made by a row of the form $(g_x,g_y)$ it would give a counterexample to the two dimensional Jacobian Conjecture. However a theorem of Magnus ([6]) shows that this is not possible. Namely his theorem asserts that if $(f,g)$ is a Keller map and either the degree of $f$ or $g$ is a prime number, then $(f,g)$ is invertible. Since $x+x^2y$ has degree $3$ it cannot be used to construct a two dimensional counterexample to the Jacobian Conjecture.
Nevertheless one could hope that $x+x^2y$ can be extended to a Keller map in dimension $3$. However at first glance I assumed that this was too much to hope for. So looking at Alpöge's example I wondered if it would be possible to construct a counterexample to the three dimensional Jacobian Conjecture, by taking as third component $F_3=x+x^2y+x^3z$, which I considered as the {\em simplest} natural generalization of the polynomial $x+x^2y$ (observe that the gradient of $x+x^2y+xz$ is not unimodular and that $x+x^2y+x^2z$ can be transformed into $x+x^2y$ by making the coordinate change $y\rightarrow y+z$).

Making every time the {\em simplest and most likely choice} will be the red thread throughout this story. The next question is: how to choose $F_1$ and $F_2$? Again guided by {\em the simplest choice} principle we choose both $F_1$ and $F_2$ to be  {\em linear in $z$}. So now the question we want to solve is:

Do there exist polynomials $F_1=a+bz, F_2=c+dz$, where $a,b,c,d$ are polynomials in $\mathbb{C}[x,y]$ such that $detJ(F_1,F_2,F_3)$ is a non-zero constant, where $F_3=x+x^2y+x^3z$? If we succeed in finding them we found the desired counterexample since $F_3=x(1+xy+x^2z)$ is clearly reducible. One more observation before we start our search: we assume that it is most likely for $b$ and $d$ to be linearly independent over $\mathbb{C}$, namely if for example $b=\rho d$, for some $\rho\in\mathbb{C}$, then replacing $F_1$ by $F_1-\rho F_2$ we get a counterexample with first component in $x$ and $y$ only, which we assume to be unlikely. Therefore from now on we will look for $b$ and $d$ which are linearly independent over $\mathbb{C}$.

\section{The coefficient of $z^2$ in the Jacobian determinant}

Let $F=(a+bz,c+dz,x+x^2y+x^3z)$, with $a,b,c,d\in\mathbb{C}[x,y]$. Then

$$
JF=\begin{bmatrix}
a_x+b_xz&a_y+b_yz&b\\
c_x+d_xz&c_y+d_yz&d\\
1+2xy+3x^2z&x^2&x^3
\end{bmatrix}
$$
The determinant of this matrix is a polynomial of degree $2$ in $z$. We first compute the coefficient of $z^2$ and investigate what is means that this coefficient equals zero. This gives the equation

$$x^3(b_xd_y-b_yd_x)-d(-3x^2b_y)+b(-3x^2d_y)=0$$

where $b_x,b_y$ etcetera denote the partial derivatives of $b$ with respect to $x$ respectively $y$. The above equation is equivalent to

$$(x\partial_x-3)(b)d_y=(x\partial_x-3)(d)b_y\,\,\,\,\,(1)$$

So we need to find solutions to this equation, not necessarily all and try them out. To get an idea how such solutions look like we start looking for solutions which are polynomials in $x$ and $y$ but with $x$-degree at most $1$, so

$$b=b_0(y)+b_1(y)x, d=d_0(y)+d_1(y)x$$

Substituting in (1) and comparing the coefficients of $x^0,x^1,x^2$ respectively gives the following equations

$$b_0{d_0}_y-{b_0}_yd_0=0\,\,\,\,(1.0)$$

$$-2b_1{d_0}_y-3b_0{d_1}_y+2d_1{b_0}_y+3d_0{b_1}_y=0\,\,\,\,(1.1)$$

$$b_1{d_1}_y-d_1{b_1}_y=0\,\,\,\,(1.2)$$

From (1.2) it follows that $(b_1/d_1)_y=0$ if $d_1\neq 0$, so $b_1$ and $d_1$ are linearly dependent, hence $b_0$ and $d_0$ cannot both be zero (otherwise $b$ and $d$ are linearly dependent). So we may assume that $b_0$ is non-zero. From (1.0) we get $(d_0/b_0)_y=0$, so $d_0=\lambda_0 b_0$ for some $\lambda_0\in\mathbb{C}$. Replacing $F_2$ by $F_2-\lambda_0 F_1$ we may assume $d_0=0$. Then equation (1.1) becomes

$$ -3b_0{d_1}_y+2d_1{b_0}_y=0\,\,\,\,$$

\noindent or equivalently $(b_0^2/d_1^3)_y=0$. So $b_0^2=\sigma d_1^3$ for some $\sigma\in\mathbb{C}$. It follows that $b_0=u^3$ and $d_1=\tau u^2$ for some $\tau\in\mathbb{C}$ and $u\in\mathbb{C}[y]$. By (1.2) $b_1=\lambda_1 d_1$, for some $\lambda_1\in\mathbb{C}$. So we get  $b=u^3+\lambda_1\tau u^2x$ and $d=\tau u^2 x$. Replacing $F_1$ by $F_1-\lambda_1 F_2$ we may assume $b=u^3$ and $d=\tau u^2 x$. Finally, replacing $F_2$ by $\tau^{-1}F_2$ we may assume $d=u^2x$.

So $b=u^3$ and $d=xu^2$ are linearly independent solutions of (1), where $u$ is an arbitrary (non-zero) polynomial in $y$. However one easily verifies that for {\em every} polynomial $u(x,y)\in\mathbb{C}[x,y]$ the polynomials $b=u(x,y)^3$ and $d=xu(x,y)^2$ are linearly independent solutions of (1) (except if $u$ is a constant multiple of $x$). Therefore from now on we assume that $b=u(x,y)^3$ and $d=xu(x,y)^2$, where we still have to choose $u(x,y)$. By this choice the $z^2$ coefficient of $detJF$ equals zero.

\section{The coefficient of $z$ in the Jacobian determinant}

The coefficient of $z$ in $detJF$ is equal to

$$x^3(a_xd_y+b_xc_y-a_yd_x-b_yc_x)-d(x^2b_x-3x^2a_y-b_y(1+2xy))$$

$$+b(x^2d_x-3x^2c_y-(1+2xy)d_y)$$

Setting this equation equal to zero and substituting $b=u^3,d=xu^2$ gives

$$x^3(a_x2xuu_y+3u_xu^2c_y-a_y(u^2+2xuu_x)-3u^2u_yc_x)$$

$$-xu^2(3x^2u^2u_x-3x^2a_y-3u^2u_y(1+2xy))$$

$$+u^3(x^2u^2+2x^3uu_x-3x^2c_y-2xuu_y(1+2xy))=0$$

We deduce that $x^2$ divides $xu^4u_y$, so $x$ divides $u^4u_y$. Now we {\em assume} that $x$ does not divide $u$. This is a choice we make. Then $x$ divides $u_y$. Write $u=u_0(y)+x\tilde{u}(x,y)$.The fact that $x$ divides $u_y$ implies that ${u_0}_y=0$. Since $x$ does not divide $u$ we have $u_0=\mu$ a non-zero constant in $\mathbb{C}$. We may assume $\mu=1$ (divide $F_1$ by $\mu^3$ and $F_2$ by $\mu^2$). So $u=1+x\tilde{u}(x,y)$. Now we make the simplest, most likely choice for $\tilde{u}(x,y)$, namely $\tilde{u}(x,y)=\lambda y$, with $\lambda$ a non-zero constant in $\mathbb{C}$. So from now on we have

$$u=1+\lambda xy, \lambda\in\mathbb{C}^*$$

So $u_x=\lambda y, u_y=\lambda x$. Substituting these values in the above equation and simplifying the whole expression and dividing by $x^2u$ finally gives

\[
2\lambda x^3 a_x
+ (2xu - 2\lambda x^2y)a_y
- 3\lambda x^2 u c_x
+ (3\lambda xyu - 3u^2)c_y
+ (\lambda xyu^3 + \lambda u^3 + u^4)
= 0.
\]

Observe that
\[
2xu - 2\lambda x^2y
= 2x(u-\lambda xy)
= 2x,
\]
\[
3\lambda xyu - 3u^2
= 3u(\lambda xy-u)
= -3u,
\]
and
\[
\lambda xyu^3+\lambda u^3+u^4
= (u-1)u^3+\lambda u^3+u^4
= (\lambda-1)u^3+2u^4.
\]

So the above equation becomes
\[
2\lambda x^3 a_x
+2x a_y
-3\lambda x^2u c_x
-3u c_y
+(\lambda-1)u^3+2u^4
=0.
\]

So
\[
2x\bigl(\lambda x^2a_x+a_y\bigr)
-3u\bigl(\lambda x^2c_x+c_y\bigr)
=
-2u^4+(1-\lambda)u^3.
\]

If we define
\[
D:=\lambda x^2\frac{\partial}{\partial x}
+\frac{\partial}{\partial y},
\]
then we get
\[
2xD(a)-3uD(c)
=
-2u^4+(1-\lambda)u^3.
\tag{2}
\]

To find solutions \(a\) and \(c\) of (2), we observe
that every polynomial of \(\mathbb{C}[x,y]\) can be written
uniquely in the form
\[
\sum_{i\geq 1}\sigma_i s_i(u)x^i
+
\sum_{j\geq 1}\tau_j t_j(u)y^j
+
w(u),
\]
where
\[
\sigma_i,\tau_j\in\mathbb{C}
\]
and \(s_i(u),t_j(u),w(u)\) are polynomials in \(u\) over
\(\mathbb{C}\).

Observe that
\[
D(u)=D(\lambda xy)=\lambda(\lambda x^2y+x)=\lambda xu.
\]

So
\[
D\bigl(s_i(u)x^i\bigr)
=
s_i'(u)\lambda xu\,x^i
+s_i(u)i\lambda x^{i+1}
\]
and hence
\[
D\bigl(s_i(u)x^i\bigr)
=
\lambda\bigl(us_i'(u)+is_i(u)\bigr)x^{i+1}.
\]

Also,
\[
D\bigl(t_j(u)y^j\bigr)
=
t_j'(u)\lambda xu\,y^j
+t_j(u)j y^{j-1}.
\]
and hence
\[
D\bigl(t_j(u)y^j\bigr)
=\bigl(t_j'(u)(u-1)u+jt_j(u)\bigr)y^{j-1}
\]
Since the right-hand side of (2) is of the form \(w(u)\),
it follows that the simplest possible solution of (2)
is of the form
\[
a=y^n p(u),
\qquad
c=y^m q(u),
\]
with \(n,m\geq 1\) and
\[
p(u),q(u)\in\mathbb{C}[u].
\]

Then
\[
2xD(a)
=
2x\left[p'(u)(u-1)u+n p(u)\right]y^{n-1}.
\]

% --- Arno-02 ---

Since
\[
xy=\frac{u-1}{\lambda},
\]
we get
\[
2xD(a)
=
\left[
2p'(u)\left(\frac{1}{\lambda}\right)(u-1)^2u
+
2n\left(\frac{1}{\lambda}\right)(u-1)p(u)
\right]y^{n-2},
\]
and
\[
3uD(c)
=
\left[
3q'(u)(u-1)u^2
+
3muq(u)
\right]y^{m-1}.
\]

Since the right-hand side of (2) does not contain terms
of the form
\[
\tau_j t_j(u)y^j,
\qquad j\geq1,
\]
we choose
\[
n=2,\qquad m=1.
\]

Hence we can write
\[
a=y^2p(u),
\qquad
c=yq(u),
\]
for some \(p(u),q(u)\in\mathbb{C}[u]\).

Then (2) becomes
\begin{align}
2p'(u)\left(\frac{1}{\lambda}\right)(u-1)^2u
&+
4\left(\frac{1}{\lambda}\right)(u-1)p(u)
-
3q'(u)(u-1)u^2
-
3u q(u) \nonumber \\
&=
-2u^4+(1-\lambda)u^3.
\tag{3}
\end{align}

It follows that \(u\) divides \(p(u)\). Since the right-hand
side of (3) has degree \(4\) in \(u\), the simplest choice for
\(p(u)\) and \(q(u)\) is that both have degree at most \(2\).
Since \(u\) divides \(p(u)\), we can write
\[
p(u)=p_2u^2+p_1u
\]
and
\[
q(u)=q_2u^2+q_1u+q_0,
\qquad p_i,q_i\in\mathbb{C}.
\]

Substituting in (3) and equating the coefficients of the
powers of \(u\) gives the following equations:
\[
-6q_2+\frac{4p_2}{\lambda}+2=0,
\]
\[
3q_2-3q_1-\frac{4p_2}{\lambda}
+\frac{2p_1}{\lambda}
+\lambda-1=0,
\]
\[
-3q_0-\frac{2p_1}{\lambda}=0.
\]

Solving \(q_2,q_1,q_0\) from these equations gives
\[
q_2
=
\frac{2}{3}\left(\frac{1}{\lambda}\right)p_2
+\frac{1}{3},
\tag{3.1}
\]
\[
q_1
=
-\frac{2}{3}\left(\frac{1}{\lambda}\right)p_2
+
\frac{2}{3}\left(\frac{1}{\lambda}\right)p_1
+
\frac{1}{3}\lambda,
\tag{3.2}
\]
\[
q_0
=
-\frac{2}{3}\left(\frac{1}{\lambda}\right)p_1.
\tag{3.3}
\]

% --- Arno-03 ---

\section*{The coefficient of \(z^0\) in the Jacobian determinant}

Finally, we have to solve for \(p_1,p_2\) and \(\lambda\)
from setting the coefficient of \(z^0\) in the Jacobian
determinant equal to a non-zero constant.

In other words, using the formulas for \(a,b,c,d\) above
and the formulas (3.1), (3.2) and (3.3), we have to solve
for \(p_1,p_2\) and \(\lambda\) from the equation
\[
\det
\begin{pmatrix}
\lambda y^3p'(u)
&
y\bigl(2p(u)+(u-1)p'(u)\bigr)
&
u^3
\\[2mm]
\lambda y^2q'(u)
&
q(u)+(u-1)q'(u)
&
xu^2
\\[2mm]
\frac{2}{\lambda}u+\left(1-\frac{2}{\lambda}\right)
&
x^2
&
x^3
\end{pmatrix}
\in\mathbb{C}^{*}.
\]

We rewrite this determinant. First multiply the first row of this matrix by$1/y$ and the third row by
$\lambda^2y^2$. Next multiply the first column by $1/\lambda^2y^2$ and the third column by $y$. These operations don't change the determinant. So we get
\[
\det
\begin{pmatrix}
\frac{p'(u)}{\lambda}&2p(u)+yu_yp'(u)&u^3\\[2mm]
\frac{q'(u)}{\lambda}&q(u)+yu_yq'(u)&yxu^2\\[2mm]
1+2xy&\lambda^2x^2y^2&\lambda^2x^3y^3
\end{pmatrix}
\]    
Now use that $u_y=\lambda x$, $\lambda xy=u-1$ and $xy=\frac{u-1}{\lambda}$.
This gives
\[
det(u):=\det
\begin{pmatrix}
\frac{p'(u)}{\lambda}&2p(u)+(u-1)p'(u)&u^3\\[2mm]
\frac{q'(u)}{\lambda}&q(u)+(u-1)q'(u)&\frac{1}{\lambda}(u-1)u^2\\[2mm]
1+2(\frac{u-1}{\lambda})&(u-1)^2&\frac{1}{\lambda}(u-1)^3
\end{pmatrix}
\]
So $det(u)$ is a polynomial of degree at most six. Now we impose
necessary conditions on $p_1,p_2$ and $\lambda$ in order for this determinant
to be a non-zero constant. Again we make the simplest choices i.e.we choose conditions which give the simplest equations. The first such a condition is that the coefficient of $u^6$
has to be zero. This gives
\[
-\frac{2p_2q_2}{\lambda^2}
-\frac{4q_2}{\lambda}
+\frac{6p_2}{\lambda^2}
=0.
\]
Substituting (3.1) and simplifying gives
\[
p_2^2-2\lambda p_2+\lambda^2=0.
\]

So
\[
p_2=\lambda,
\]
and hence by (3.1),
\[
q_2=1.
\]

Next, setting the coefficient of $u$ equal to zero gives

\[
\frac{p_1q_1}{\lambda^2}-\frac{2p_2q_0}{\lambda^2}+\frac{3p_1q_0}{\lambda^2}=0
\]

Using (3.2), (3.3) and $p_2=\lambda$ we obtain
\[
p_1=\frac{\lambda^2}{4}+\frac{\lambda}{2}
\]

Substituting these values of $p_1$ and $p_2$ in (3.2) respectively (3.3) gives

$$q_1=\frac{\lambda}{2}-\frac{1}{3},  \,\,\, q_0=-\frac{\lambda}{6}-\frac{1}{3}$$

The coefficient of $u^0$ equals $\frac{-p_1 q_0}{\lambda^2}$. Substituting
the values of $p_1$ and $q_0$ gives the constant
$$\frac{\lambda}{24}+\frac{1}{6}+\frac{1}{6\lambda}$$
Now observe that another necessary condition for the determinant to be constant is 
that we get the same value if we substitute $u=1$ in the determinant. The result is
$-q(1)=-(q_0+q_1+q_2)=-(\frac{\lambda}{3}+\frac{1}{3})$
So equating these two constants gives
$$\frac{\lambda}{24}+\frac{1}{6}+\frac{1}{6\lambda}=-(\frac{\lambda}{3}+\frac{1}{3})$$
or equivalently $9\lambda^2+12\lambda+4=0$, which gives
$$\lambda=-\frac{2}{3}$$

So we get

\[
p(u)
=
-\frac{2}{3}u^2-\frac{2}{9}u,
\]
\[
q(u)
=
u^2-\frac{2}{3}u-\frac{2}{9},
\]
and
\[
u=1-\frac{2}{3}xy.
\]

With these values of $\lambda, p(u)$ and $q(u)$ one easily verifies that $det(u)$ is a non-zero constant, more precisely
it is equal to $-\frac{1}{9}$.
So we found a counterexample to the Jacobian Conjecture!

Summarizing, we get
\[
F_1=u^3z+y^2p(u),
\qquad
F_2=xu^2z+yq(u),
\qquad
F_3=x+x^2y+x^3z
\]
with $u,p(u),q(u)$ as above.

\section{Some final remarks}

The map we constructed above is essentially the map found by Alp\"oge. More precisely, 
if we denote by
\[
A=(A_1,A_2,A_3)
\]
Alp\"oge's counterexample,
and let
\[
L_1=
\left(
\frac{1}{2}x,
-\frac{4}{3}y,
-8z
\right),
\qquad
L_2=
\left(
-\frac{1}{8}x,
-\frac{1}{12}y,
z
\right),
\]
then
\[
L_2\circ A\circ L_1=F.
\]

The second remark concerns the polynomial $x+x^2y$. Our construction started with
this polynomial and the idea that I wondered if it would be possible to construct
a counterexample to the three dimensional Jacobian Conjecture by taking as its third 
component the polynomial $x+x^2y+x^3z$ which I considered as the simplest natural
generalization of the polynomial $x+x^2y$. I also assumed that it was too much to hope
for that $x+x^2y$ could be extended to a three dimensional counterexample. However in
retrospect we can now conclude that $x+x^2y$ {\em can} be extended to such an example.
Namely, if we take 
$$E=(x,y-xz,z)$$
then the third component of the Keller map $F\circ E$ is equal to $x+x^2y$.

\section{Acknowledgements}

I want to thank Jan Schoone for sending me the counterexample and for checking a first
version of this paper. Also I want to thank Stefan Maubach for setting me on the right track.
Then as always Michiel de Bondt for his help with Maxima and his many very useful comments. Engelbert Hubbers for making several corrections in the text and last, but certainly not least,
Joris Vergeest without whose help this note was probably never published.

\section{References}

[1] L. Alp\"oge, A counterexample to the Jacobian Conjecture in dimension three, announcement July 19, 2026.

\noindent [2] A. van den Essen, Polynomial Automorphisms and the Jacobian Conjecture, Progress in Mathematics, 190, Birkh\"auser, Basel 2000.

\noindent [3] A. Gallagher, An infinite family of counterexamples to the Jacobian Conjecture in dimension three: every generic fiber of degree $n\geq 3$ occurs, preprint July 20, 2026, Zenodo.

\noindent [4] S. Gao, Counterexamples to the Jacobian Conjecture in dimensions greater than two,
arXiv: 2608.00222v1

\noindent [5] O. Keller, Ganze Cremona Transformationen, Monatsh. Math. Phys. 47 (1939), 299-306.

\noindent [6] A. Magnus, On polynomial solutions to a differential equation, Math. Scand., 3 (1955), 255-260

\noindent [7] D. Speyer, The geomretry and structure of Gallagher's counterexamples to the Jacobian Conjecture, July 23, 2026.

\noindent [8] T. Tao, A digestion of the Jacobian Conjecture counterexample, blog post, July 21, 2026

\noindent essen@math.ru.nl, arno.vd.essen@gmail.com

\end{document}